\documentclass{amsart}
\usepackage{amsmath,amsthm,amsfonts}
\usepackage{csquotes}

\newcommand{\coh}[1]{\mathsf{Coh}({#1})}
\newcommand{\ddp}[2]{\partial_{#1}^{[#2]}}
\newcommand{\hol}[1]{\mathsf{Hol}({#1})}
\newcommand{\mcat}[1]{\mathsf{Mod}({#1})}
\newcommand{\n}{\mathbb{N}_0}
\newcommand{\q}{\mathbb{Q}}
\newcommand{\z}{\mathbb{Z}}
\newcommand{\mv}[1]{\underline{#1}}
\newcommand{\aff}{\mathbb A}
\DeclareMathOperator{\coker}{coker}
\DeclareMathOperator{\der}{Der}
\DeclareMathOperator{\enm}{End}
\DeclareMathOperator{\ext}{Ext}
\DeclareMathOperator{\fchar}{char}

\DeclareMathOperator{\hmm}{Hom}

\theoremstyle{plain}
\newtheorem{thm}{Theorem}
\newtheorem{lem}[thm]{Lemma}
\newtheorem{prop}[thm]{Proposition}

\begin{document}

\title{Coherent rings of differential operators}
\author{Eivind Eriksen}
\address{BI Norwegian Business School, Department of Economics, N-0442
Oslo, Norway}
\email{eivind.eriksen@bi.no}
\date{\today}

\begin{abstract}
We consider the following question: When are rings of differential operators 
coherent? If $A$ is a finitely generated smooth domain over a field $k$ of 
characteristic $0$, then the ring $D$ of differential operators on $A$ is a 
Noetherian ring and a finitely generated $k$-algebra. However, when $k$ has 
characteristic $p > 0$ or when $A$ is singular, this is no longer true. In 
fact, Bernstein, Gelfand and Gelfand showed that for the cubic cone $A = 
k[x,y,z]/(x^3 + y^3 + z^3)$, the ring $D$ is neither Noetherian nor finitely 
generated if $\fchar(k) = 0$, and the same is true for the polynomial ring 
$A = k[x_1, \dots, x_n]$ if $\fchar(k) = p > 0$. In this paper, we prove that 
the ring $D$ of differential operators on a finitely generated, smooth and 
connected algebra $A$ over a field $k$ of characteristic $p > 0$ is coherent, 
and conjecture that same holds for the cubic cone in characteristic $0$. We 
argue that the question of coherence is the more fundamental one, and use
some interesting results of Bavula to study holonomic $D$-modules on $A = 
k[x_1, \dots, x_n]$ in characteristic $p > 0$. 
\end{abstract}


\maketitle

\section{Introduction}

Let $A$ be a finitely generated commutative algebra over a field $k$, and
consider the ring $D$ of $k$-linear differential operators on $A$ in the
sense of Grothendieck \cite{grot67}. Any system of linear differential
equations defined over $A$ has the form
\begin{align*}
P_{11}(u_1) &+ P_{12}(u_2) + \dots + P_{1n}(u_n) = 0 \\
P_{21}(u_1) &+ P_{22}(u_2) + \dots + P_{2n}(u_n) = 0 \\
\vdots & \\
P_{m1}(u_1) &+ P_{m2}(u_2) + \dots + P_{mn}(u_n) = 0
\end{align*}
where $P_{ij} \in D$ are differential operators for $1 \le i \le m, \; 1 \le
j \le n$ and $u_1, u_2, \dots, u_n$ are unknown functions. To this system,
we can associate the finitely presented left $D$-module $M = \coker(\phi)$,
given by
    \[ 0 \gets M \gets D^n \xleftarrow{\phi} D^m \]
where $\phi: D^m \to D^n$ is given by right multiplication by the $m \times
n$ matrix $( P_{ij} )$ with coefficients in $D$. It is therefore reasonable 
to define an algebraic $D$-module to be a finitely presented left $D$-module 
$M$. With this definition, we may identify $\hmm_D(M,\mathcal S)$ with the 
set of solutions of the above system of differential equations with values
in a left $D$-module $\mathcal S$.

If $A$ is a smooth integral domain over a field $k$ of characterstic $0$,
then $D$ is a simple Noetherian ring, and a finitely generated $k$-algebra.
In this case, any finitely generated left $D$-module is finitely presented.
However, $D$ is neither finitely generated nor Noetherian in general; see
for instance Smith \cite{smit87}. In fact, when $A$ is either of the rings 
\begin{enumerate}
\item $A = k[x,y,z]/(x^3+y^3+z^3)$ when $k$ has characteristic $0$
\item $A = k[x_1, \dots, x_n]$ when $k$ has characteristic $p > 0$
\end{enumerate}
then the ring $D$ of differential operators on $A$ is not a Noetherian ring, 
and not a finitely generated $k$-algebra. The study of algebraic $D$-modules 
are in these cases considered hopeless, but this is not necessarily so; if 
$D$ is a coherent ring, then the category of finitely presented left 
$D$-modules, or coherent left $D$-modules, have good properties.

The main result in this paper is that the ring $D$ of differential operators 
on a finitely generated, smooth and connected commutative algebra $A$ over a 
field $k$ of characteristic $p > 0$ is coherent. We conjecture that this is
also the case for the ring of differential operators on the cubic cone $A = 
k[x,y,z]/(x^3 + y^3 + z^3)$ when $k$ has characteristic $0$. However, this 
is an open question, as far as we know.

In Bavula \cite{bavu09}, the author studies finitely generated and finitely
presented left $D$-modules when $D$ is the ring of differential operators on
the polynomial ring $A = k[x_1, \dots, x_n]$ over a field $k$ of 
characteristic $p > 0$. His results show that finitely presented $D$-modules 
give a far more reasonable theory than finitely generated $D$-modules. Our 
results show that $D$ is a coherent ring in this case, and a $D$-module is 
therefore coherent if and only if it is finitely presented. Bavula's results 
fit very nicely with our point of view, that coherent modules is the 
\enquote{correct} notion for $D$-modules. 

The results in Bavula \cite{bavu09} give a classification of the holonomic 
$D$-modules when $D$ is the ring of differential operators on $A = k[x_1, 
\dots, x_n]$ in characteristic $p > 0$. The result is that there are very 
few holonomic $D$-modules, and they are trivial as differential operators 
since they are given by multiplicative operators of order $0$. As a 
comment to these results, we show that when $p = 2$ and $n = 1$, then the 
coherent left $D$-module $M = D/D \cdot \partial$ is not holonomic, but 
has dimension $\dim(M) = 2$ and multiplicity $e(M) = 1/2$.

\section{Coherent rings and modules}

Let $R$ be an associative ring. A left $R$-module $M$ is \emph{coherent} if
it is a finitely generated $R$-module with the property that any finitely
generated submodule of $M$ is finitely presented. We write $\coh R \subseteq
\mcat R$ for the full subcategory of coherent left $R$-modules.

We recall some fundamental results for coherent modules. Proofs of these
results are given in Chaper 2 of Glaz \cite{glaz89} in the commutative case
(the same proofs hold when $R$ is any associative ring); see also Exercise
I.2.11-12 in Bourbaki \cite{bour-ac61i}.

\begin{lem} \label{l:cohproj}
Any finitely generated projective left $R$-module is finitely presented.
\end{lem}
\begin{proof}
This follows from the proof of Theorem 2.1.4 in Glaz \cite{glaz89}.
\end{proof}

\begin{lem} \label{l:cohmod}
Any finitely generated submodule of a coherent module is coherent, and the
full subcategory $\coh R \subseteq \mcat R$ of coherent modules is an exact
Abelian subcategory that is closed under extensions.
\end{lem}
\begin{proof}
This follows from the proof of Theorem 2.2.1, Corollary 2.2.2 and Corollary
2.2.3 in Glaz \cite{glaz89}, since a full subcategory of an Abelian category
is an exact Abelian subcategory if and only if it is closed under kernels,
cokernels and finite direct sums.
\end{proof}

We say that $R$ is a \emph{left coherent ring} if $R$ is a coherent as a left
$R$-module, or equivalently, if any finitely generated left ideal in $R$ is
finitely presented. If follows that a left Noetherian ring is left coherent.
We recall that $R$ is \emph{left semi-hereditary} if any finitely generated
left ideal in $R$ is projective. It follows from Lemma \ref{l:cohproj} that
a left semi-hereditary ring is left coherent.

\begin{prop} \label{p:cohflatlim}
Let $I$ be a directed partially ordered set, let $(R_i)_{i \in I}$ be a
direct system of associative rings, and let
    \[ R = \varinjlim \, R_i \]
be its direct limit. If $R_i$ is a left coherent ring for all $i \in I$ and
$R_i$ is a flat right $R_j$-module for all $i \ge j$ in $I$, then $R$ is a
left coherent ring.
\end{prop}
\begin{proof}
Let $a \subseteq R$ be a finitely generated left ideal. Then there exists an
index $j$ and a finitely generated left ideal $a_j \subseteq R_j$ such that
$R \otimes_{R_j} a_j \cong a$, and by the left coherence of $R_j$, it follows
that $a_j$ is finitely presented. We choose a finite presentation $D_j^m \to
D_j^n \to a_j \to 0$, and consider the sequence
    \[ R \otimes_{R_j} R_j^m \to R \otimes_{R_j} R_j^n \to R \otimes_{R_j}
    a_j \to 0 \]
of left $R$-modules. Since $R_i$ is a flat right $R_j$-module for all $i \ge
j$, it follows that $R$ is a flat right $R_j$-module, and therefore this
gives a finite presentation of the left $R$-module $a$.
\end{proof}

\begin{lem}
If $R$ is left coherent, then a left $R$-module $M$ is coherent if and only
if it is finitely presented. In this case, there is a free resolution
    \[ 0 \gets M \gets L_0 \gets L_1 \gets \dots \gets L_i \gets \dots \]
of $M$, where $L_i$ is free of finite rank for all $i \ge 0$.
\end{lem}
\begin{proof}
If $M$ is finitely presented, then $M \cong \coker(\phi)$ for a morphism
$\phi: R^n \to R^m$ of left $R$-modules. Since $R$ is coherent, the same
holds for $R^n$ and $R^m$, and $M$ is coherent by Lemma \ref{l:cohmod}.
Conversely, if $M$ is a left coherent $R$-module, then $M$ is finitely
generated, and there is an exact sequence $0 \to N \to R^m \to M \to 0$.
Since $R^n$ and $M$ are coherent, the same holds for $N = \ker(R^m \to M)$
by Lemma \ref{l:cohmod}. In particular, $N$ is finitely generated and $M$
is finitely presented. By the coherence of $N$, it also follows that $N$
is finitely presented, and an inductive argument shows that we can extend
the finite presentation of $M$ to a free resolution
    \[ 0 \gets M \gets L_0 \gets L_1 \gets \dots \gets L_i \gets \dots \]
of $M$, with $L_i$ free of finite rank for all $i \ge 0$.
\end{proof}

The notion of a right coherent module and of a right coherent ring can be
defined similarly, and by symmetry, the results in this section also hold
for right modules. We say that $R$ is a \emph{coherent ring} if it is left
and right coherent. The polynomial ring $R = k[x_1, x_2, \dots ]$ in an
infinite number of variables $x_1, x_2, \dots$ over a field $k$ is an
example of a coherent ring that is not Noetherian.

\section{The ring of differential operators on a polynomial ring}

Let $A = k[x_1, x_2, \dots, x_n]$ be the polynomial ring in $n$ variables
over a field $k$. We consider the ring $D = D(A)$ of $k$-linear differential
operators on $A$, in the sense of Grothendieck \cite{grot67}. This is a
filtered ring, equipped with the order filtration
    \[ A = D^0 \subseteq D^1 \subseteq \dots \subseteq D^i \subseteq \dots
    \subseteq D \quad \text{with} \quad D = \bigcup_{i \ge 0} \, D^i \]
Let us write $\mu: A \otimes_k A \to A$ for the multiplication map, given by
$\mu(a \otimes b) = ab$, and $J = \ker(\mu)$ for its kernel, which acts on
$\enm_k(A)$ in the natural way. The the set of differential operators of 
order at most $i$ is given by $D^i = \{ P \in \enm_k(A): J^i \cdot P = 0 
\}$. 

Let us describe the ring $D$ in concrete terms. We consider the partial 
derivations $\partial_i = \partial / \partial x_i \in \der_k(A)$ for $1 \le i 
\le n$, and define their \emph{divided powers} $\ddp{i}{r}: A \to A$ to be 
the $k$-linear operators given by
    \[ \ddp{i}{r}(\mathbf x^m) = \binom{m_i}{r} \, \mathbf x^{m - r
    \epsilon_i} \]
for all multi-indices $m = (m_1, m_2, \dots, m_n) \in \n^n$ and for all
integers $r \ge 0$, where we use multi-index notation
    \[ \mathbf x^m = x_1^{m_1} x_2^{m_2} \cdots x_n^{m_n} \]
and write $m - r \epsilon_i = (m_1, \dots, m_i-r, \dots, m_n)$. Notice that
the binomial coefficients in $k$ are the canonical images of the usual
integer-valued binomial coefficients. The name divided powers come from the
fact that $r! \, \ddp{i}{r} = \partial_i^r$. The following result is
well-known, see for instance Section 4 in Bavula \cite{bavu09}:

\begin{lem} \label{l:weylp}
The ring $D = D(k[x_1,\dots,x_n])$ is the subalgebra of $\enm_k(A)$ generated
by $x_i$ and the divided powers $\ddp{i}{r}$ for $1 \le i \le n$ and $r \ge
1$. These generators have relations given by
    \[ \left[ x_i, x_j \right] = \left[ \ddp{i}{r}, \ddp{j}{s} \right] = 0,
    \quad \ddp{i}{r} \, \ddp{i}{s} = \binom{r+s}{r} \ddp{i}{r+s}, \quad
    \left[ \ddp{i}{r}, x_j \right] = \delta_{ij} \, \ddp{i}{r-1} \]
for all $1 \le i,j \le n$ and all $r,s \ge 1$.
\end{lem}

If $\fchar(k) = 0$, then $D = A_n(k)$ is the $n$'th Weyl algebra, which is a
simple Noetherian ring, generated by $\{ x_1, x_2, \dots, x_n, \partial_1,
\dots, \partial_n \}$. If $\fchar(k) = p > 0$, then it is known that $D$ is
not Noetherian and not a finitely generated $k$-algebra. We claim that $D$ 
is a coherent ring. In fact, we shall prove a more general result in the next 
section.

In any characteristic, we have a finite dimensional filtration $\{ B^i \}$
of the ring $D$, given by the $k$-linear spaces
    \[ B^i = \left\{ \sum_{|m|+ |r| \le i} c_{m,r} \, \mathbf x^m \,
    \mathbf \partial^{[r]}: c_{m,r} \in k \text{ for all } m,r \in \n^n
    \right\} \]
for $i \ge 0$, where $\mathbf \partial^{[r]} = \partial_1^{[r_1]} \cdots
\partial_n^{[r_n]}$. We follow Bavula \cite{bavu09} and call this filtration
the \emph{canonical filtration} of $D$. When $\fchar(k) = 0$, it coincides 
with the usual Bernstein filtration. Notice that $\dim_k B^i < \infty$ for 
all $i \ge 0$, since we have
    \[ \dim_k B^i/B^{i-1} = \binom{2n+i-1}{i} \]
Moreover, we have that $B_i = 0$ for $i < 0$, that $B^i \subseteq B^{i+1}$
and $B^i \cdot B^j \subseteq B^{i+j}$ for all $i,j \ge 0$, and that $\cup_i
\, B^i = D$.

\section{Coherent rings of differential operators}

Let $k$ be a field of characteristic $p > 0$, and let $A$ be a finitely
generated, smooth and connected commutative algebra over $k$. We shall use
the following construction, introduced in Section 3 of Chase \cite{chas74}:
Let $A_r \subseteq A$ be the $k$-subalgebra generated by $\{ a^{p^r}: a \in
A \}$ for all $r \ge 0$, and consider the chain
    \[ A = A_0 \supseteq A_1 \supseteq \dots \supseteq A_r \supseteq
    \dots \]
of $k$-algebras. We define $D_r = \enm_{A_r}(A) \subseteq \enm_k(A)$ for
$r \ge 0$, and identify $A$ with $D_0 = \enm_A(A)$. From Lemma 3.3 in Chase
\cite{chas74}, it follows that
    \[ A = D_0 \subseteq D_1 \subseteq \dots \subseteq D_r \subseteq
    \dots \subseteq D \quad \text{and} \quad \bigcup_{r \ge 0} \, D_r = D \]
where $D = D(A)$ is the ring of $k$-linear differential operators on $A$.

\begin{prop} \label{p:proj-charp}
Let $k$ be a field of characteristic $p > 0$, and let $A$ be a finitely
generated, smooth and connected commutative algebra over $k$. Then we have:
\begin{enumerate}
\item $A$ is a finitely generated projective $A_r$-module for all $r \ge 0$
\item $D_r$ is Morita equivalent to $A_r$ for all $r \ge 0$
\item $D_s$ is a projective right $D_r$-module for all $r \le s$
\end{enumerate}
\end{prop}
\begin{proof}
The first part follows from Lemma 3.2 in Chase \cite{chas74}. Since $A$ is
clearly a faithful $A_r$-module, and $A_r$ is commutative, it follows that
$A$ is a progenerator, and $D_r = \enm_{A_r}(A)$ is Morita equivalent to
$A_r$. The last part follows from the proof of Proposition 3.2 in Smith
\cite{smit87}.
\end{proof}

\begin{thm} \label{t:diffcoh}
Let $A$ be a finitely generated, smooth and connected commutative algebra
over a field $k$ of characteristic $p > 0$. Then the ring $D = D(A)$ of
$k$-linear differential operators on $A$ is a coherent ring.  
\end{thm}
\begin{proof}
It follows from the comments above that the ring $D = D(A)$ of differential
operators on $A$ is the direct limit
    \[ D(A) = \varinjlim \, D_r \]
of the directed system $D_0 \subseteq D_1 \subseteq \dots$ of subrings.
Moreover, $D_r$ is Noetherian, and therefore coherent, since it is Morita
equivalent to $A_r$ by Proposition \ref{p:proj-charp}, and $D_s$ if a flat
right $D_r$-module for all $r \le s$ since it is projective by Proposition
\ref{p:proj-charp}. Hence $D = D(A)$ is a coherent ring by Proposition
\ref{p:cohflatlim}.
\end{proof}

When $A = \mathcal O(X)$ is the coordinate ring of a non-singular affine
algebraic variety $X$ of dimension $d$ over an algebraically closed field
$k$ of characteristic $p > 0$, then $D = D(A)$ has global homological
dimension $d$ by Theorem 3.7 in Smith \cite{smit87}. In case $d \le 1$, it
is therefore known that the ring $D = D(A)$ of differential operators is a
coherent ring; by definition, any hereditary ring is semi-hereditary, and
therefore coherent. As far as we know, this result is new for $d \ge 2$.

\section{Holonomic $D$-modules}

In this section, we assume that $A = k[x_1, \dots, x_n]$ is a polynomial
algebra over a field $k$ of characteristic $p > 0$, and that $D = D(A)$
is the ring of differential operators on $A$. Then it follows from Theorem
\ref{t:diffcoh} that $D$ is a coherent ring, and we consider the category
$\coh D$ of coherent left $D$-modules. 

Let $M$ be a left $D$-module. If $M$ is finitely generated, then there is
a finite dimensional $k$-linear subspace $M_0 \subseteq M$ such that $M =
D \cdot M_0$. We consider the finite dimensional filtration $\{ M_i \}$ of
$M$ given by
    \[ M_i = B^i \cdot M_0 \]
for $i \ge 0$, and define $\dim(M)$ to be the growth of the function $i
\mapsto \dim_k M_i$. We recall that for a function $f: \n \to \n$, the
growth $\gamma(f)$ is defined as
    \[ \gamma(f) = \inf \{ m: f(i) \le i^m \text{ for } m \gg 0 \} \]
This definition of $\dim(M)$ appears in Bavula \cite{bavu09}. It does not
depend on the choice of generating set $M_0$, but may depend on the choice
of finite dimensional filtration $\{ B^i \}$ of $D$. We shall therefore 
fix the canonical filtration of $D$. The definition of $\dim(M)$ given by 
Bavula resembles the definition of Gelfand-Kirillov dimension over finite 
dimensional algebras, but is better suited for coherent $D$-modules in 
positive characteristic. 

\begin{prop}[Bavula] \label{p:bern-ineq}
If $M$ is a non-zero, finitely generated left $D$-module, then $\dim(M)
\ge n$.
\end{prop}
\begin{proof}
See Theorem 4.3 in Bavula \cite{bavu09}.
\end{proof}

Let $M$ be a coherent left $D$-module. We say that $M$ is \emph{holonomic}
if $\dim(M) = n$, and define the category $\hol D$ of holonomic $D$-modules
to be the full subcategory of $\coh D$ consisting of holonomic $D$-modules.
This definition is different than the one used by Bavula, since he does not
require that $M$ is coherent. However, the definitions coincide for coherent 
$D$-modules.

To require that holonomic modules are coherent, as we do, have consequences.
This means that $A$ considered as a left $D$-module is not holonomic. For 
instance, if $A = k[x]$, then the left $D$-module $A$ can be written as $A = 
D/I$, where $I$ is the left ideal 
    \[ I = D (\partial, \partial^{[2]}, \partial^{[3]}, \dots ) \] 
Since $I$ is not finitely generated, $A = D/I$ is not finitely presented and
therefore not coherent. 

\begin{prop}
Any finitely generated submodule of a holonomic $D$-module is holonomic, and
$\hol D \subseteq \coh D$ is an exact Abelian subcategory which is closed
under extensions.  
\end{prop}
\begin{proof}
If $M$ is a holonomic $D$-module, and $N \subseteq M$ is a finitely generated
submodule, then $N$ is coherent by Lemma \ref{l:cohmod}. Moreover, if $N_0$
is a finitely dimensional subset $N_0 \subseteq N$ such that $D \cdot N_0 =
N$, then there is a finite dimensional subset $M_0 \subseteq M$ with $D \cdot
M_0 = M$ containing $N_0$, and this implies that $\dim(N) \le \dim(M) = n$.
Hence, $\dim(N) = n$ by Proposition \ref{p:bern-ineq}, and $N$ is holonomic.
If follows from Theorem 5.10 in Bavula \cite{bavu09} that $\hol D \subseteq
\coh D$ is an exact Abelian subcategory closed under extensions. 
\end{proof}

Let $M$ be a coherent left $D$-module, and let $\{ M_i \}$ be a finite 
dimensional filtration given by $M_i = B^i \cdot M_0$, where $M_0 \subseteq 
M$ is a finite dimensional linear subspace with $D \cdot M_0 = M$. Then there
is an integer $k \ge 0$ and polynomials $p_i(t) \in \q[t]$ for $0 \le i < 
p^k$ such that $\dim_k \left( M_{m p^k + i} \right) = p_i(m)$ for all $m \gg 
0$. Moreover, the polynomials $p_i(t)$ all have the same leading term, and 
are given by  
    \[ p_i(t) = \frac{e}{d!} \cdot t^d + \text{ terms of lower degree } \] 
where $d = \dim(M)$ and $e = e(M) > 0$ is the \emph{multiplicity} of $M$; see 
Theorem 5.5 of Bavula \cite{bavu09}. The function $m p^k + i \mapsto p_i(m)$ 
for $0 \le i < p^k$ and $m \ge 0$ is 
called a quasi-polynomial of period $p^k$. In general, we have that $p^{kn} 
e(M) \in \z$. If $M$ is holonomic, then $e(M)$ is a positive integer by 
Theorem 8.7 in Bavula \cite{bavu09}. By Theorem 9.6 and Corollary 6.8, we
have the following results: 

\begin{prop}
The category $\hol D$ of holonomic modules is a length category, and any 
submodule or factor module of a holonomic module is holonomic. The simple 
objects of $\hol D$ are simple considered as left $D$-modules.
\end{prop}

\begin{thm}[Bavula]
Let $A = k[x_1, \dots, x_n]$ be a polynomial ring over a field $k$ of 
characteristic $p > 0$, and let $D = D(A)$ be its ring of differential 
operators. If $k$ is algebraically closed, then the simple objects of $\hol 
D$ are given by
    \[ M(\mv \alpha) = D \otimes_A A/(x_1 - \alpha_1, \dots, x_n - \alpha_n)
    \cong D/D(x_1 - \alpha_1, \dots, x_n - \alpha_n) \]
for $\mv \alpha \in \aff^n_k$.
\end{thm}

This means that the simple holonomic $D$-modules are given by multiplication 
operators of order zero, and therefore trivial as systems of differential 
equations. Iterated extensions of these simple modules do not give more 
interesting holonomic modules. For example, when $n = 1$, we have that 
$\ext^1_D(M(\mv \alpha),M(\mv \beta)) = 0$ for all $\mv \alpha, \mv \beta \in 
\aff^1_k$. 

Hence, any linear system of \enquote{interesting} differential equations, 
given by a matrix $( P_{ij} )$ of differential operators as explained in the 
introduction, corresponds to a coherent left $D$-module that is not holonomic. 
The following example is instructive: Let $n = 1$, and let $M = D/D \cdot 
\partial$. Then there is an exact sequence 
    \[ D \xrightarrow{\cdot \partial} D \to M \to 0 \] 
In characteristic $p > 0$, the kernel of the map $D \xrightarrow{\cdot 
\partial} D$ is non-zero. For instance, if the characteristic $p = 2$, then 
the kernel is $D \cdot \partial$ since we have that 
    \[ \partial^{[m]} \cdot \partial = \begin{cases} 0, & m \text{ is odd} \\
    \partial^{[m+1]}, & m \text{ is even} \end{cases} \] 
Therefore, there is a short exact sequence of coherent left $D$-modules 
    \[ 0 \to M \to D \to M \to 0 \] 
Since $\dim(D) = 2$, we must have $\dim(M) = 2$, and $M$ is not holonomic. In 
fact, $\dim(M) = 2$ implies that the multiplicity $e(M) = e(D)/2 = 1/2$.

\bibliographystyle{amsplain}
\bibliography{eeriksen}

\end{document}